\documentclass[11pt]{article}
\usepackage{amsfonts,amssymb,amsmath,mathrsfs,bbm}
\usepackage{color,latexsym,amsfonts}

\numberwithin{equation}{section}
\newtheorem{theorem}{Theorem}[section]
\newtheorem{lemma}{Lemma}

\newtheorem{corollary}[theorem]{Corollary}

\def\<{\langle}
\def\>{\rangle}

\begin{document}

\bigskip \bigskip \noindent {\Large \textbf{Reduced critical Bellman-Harris
branching processes for small populations }}\footnote{%
\noindent The work of of V.Vatutin is supported by the Russian Science
Foundation under the grant 14-50-00005 and performed in Steklov Mathematical
Institute of Russian Academy of Sciences, Wenming Hong and Yao Ji were supported by the Natural
Science Foundation of China under the grant 11531001 and 11626245}

\noindent {
Wenming Hong\footnote{%
School of Mathematical Sciences \& Laboratory of Mathematics and Complex
Systems, Beijing Normal University, Beijing 100875, P.R. China. Email:
wmhong@bnu.edu.cn} ~ Yao Ji\footnote{%
School of Mathematical Sciences \& Laboratory of Mathematics and Complex
Systems, Beijing Normal University, Beijing 100875, P.R. China. Email:
happyjiyao@mail.bnu.edu.cn}} ~ Vladimir Vatutin \footnote{%
Steklov Mathematical Institute of Russian Academy of Sciences, 8 Gubkina
St., Moscow 119991 Russia Email: vatutin@mi.ras.ru}

\begin{center}
\begin{minipage}{12cm}
\begin{center}\textbf{Abstract}\end{center}
\footnotesize
Let $\left\{ Z(t), t\geq 0\right\} $ be a critical Bellman-Harris branching
process with finite variance for the offspring size of particles. Assuming that $0<Z(t)\leq
\varphi (t)$, where either  $\varphi
(t)=o(t)$ as $t\rightarrow \infty $ or $\varphi
(t)=at,\, a>0$, we study the structure of the process $%
\left\{ Z(s,t),0\leq s\leq t\right\} ,$ where $Z(s,t)$ is the number of
particles in the process at moment $s$  in the initial
process which either survive up to moment $t$ or have a positive
offspring number at this moment.

\bigskip

\mbox{}\textbf{Keywords:}\quad Bellman-Harris branching process, reduced process, conditional limit theorem; \\
\mbox{}\textbf{Mathematics Subject Classification}:  Primary 60J80;
secondary 60G50.

\end{minipage}
\end{center}

\section{ Introduction and main results\label{Introd}}

Let $\left\{ Z(t),t\geq 0\right\} $ be a Bellman-Harris branching process
with $Z(0)=1$ specified by the probability generating function
\begin{equation}
f(s)=\mathbf{E}s^{\xi }=\sum_{k=0}^{\infty }f_{k}s^{k}  \label{DefGener}
\end{equation}%
and the distribution $G(t)=\mathbf{P}(\tau \leq t)$ of the life-length $\tau
$ of a particle.

Introduce the following hypothesis:

\textbf{Condition A1} (Criticality)
\begin{equation}
\mathbf{E}\xi =1,\quad \sigma ^{2}:=Var\xi \in \left( 0,\infty \right) .
\notag  \label{BasicCond}
\end{equation}

\textbf{Condition A2}. The support of the distribution $G(t)$ is contained
on the integer lattice $t=0,1,2,...$ with maximal step 1 and is not
degenerate.

Let $\mu :=\mathbf{E}\tau $ and
\begin{equation*}
F(t;s):=\mathbf{E}\left[ s^{Z(t)}|Z(0)=1\right]
\end{equation*}%
be the probability generating function for the number of particles in the
process at moment $t$. It is known (see, for instance, \cite{Gol1971}) that
if Condition A1 is valid and $t^{2}\left( 1-G(t)\right) \rightarrow ~0$ as $%
t\rightarrow \infty $ then
\begin{equation}
Q(t):=1-F(t;0)=\mathbf{P}\left( Z(t)>0\right) \sim \frac{2\mu }{\sigma ^{2}t}%
\text{ \ \ as \ \ }t\rightarrow \infty  \label{SurvivalProbab1}
\end{equation}%
and, for any $\lambda \geq 0$%
\begin{equation}
\lim_{t\rightarrow \infty }\mathbf{E}\left[ e^{-2\mu \lambda Z(t)/\sigma
^{2}t}|Z(t)>0\right] =\frac{1}{1+\lambda }  \label{Yaglom0}
\end{equation}%
meaning that the limiting distribution of the scaled process $2\mu
Z(t)/\sigma ^{2}$ given $\left\{ Z(t)>0\right\} $ is exponential with
parameter 1.

In this note we study the asymptotic properties of the so-called reduced
critical Bellman-Harris process $\left\{ Z(s,t),0\leq s\leq t\right\} ,$
where $Z(s,t)$ is the number of particles at moment $s$ in the initial
process which either survive up to moment $t$ or have a positive offspring
number at this moment.

Note that reduced processes for ordinary Galton--Watson branching processes
(i.e, for the case $\mathbf{P}\left( \tau =1\right) =1)$ were introduced by
Fleischmann and Prehn \cite{FP}. Various properties of such processes were
analyzed in \cite{Athr2012},\cite{Ath2012b},\cite{Dur78},\cite{FZ}, \cite%
{HJR17}, \cite{Sam17}, \cite{Lam2003},\cite{Le2014}, \cite{Zub} and some
other papers.

Reduced critical Bellman-Harris processes were investigated by Vatutin \cite%
{Vat1979} for the single-type case and by Sagitov \cite{Sag1986} for
multitype setting.

All these papers do not consider the situation when the size of the
population at moment $n$ is bounded from above. Recently, Liu and Vatutin
\cite{LV2018} study the structure of the Galton-Watson critical reduced
process under the condition that the size of the population is bounded and
positive at the moment of observation. In the present paper we consider a
similar problem for the critical Bellman-Harris processes.

Introduce the event%
\begin{equation*}
\mathcal{H}(t):=\left\{ 0<Z(t)\leq B\varphi (t)\right\}
\end{equation*}%
where%
\begin{equation*}
B=\frac{\sigma ^{2}}{2\mu }.
\end{equation*}

Our main results are contained in two theorems which we formulate below.

\begin{theorem}
\label{T_main}Let Conditions A1-A2 be valid,
\begin{equation*}
\ \mathbf{E}\xi ^{2}\log (\xi +1)<\infty ,\ \mathbf{E}\tau ^{3}<\infty ,
\end{equation*}%
and $\varphi (t), t>0,$ be a monotone increasing function,  $\varphi (t)=o(t) $ as $t\rightarrow \infty. $ 
 If, in addition,
\begin{equation}
\lim_{t\rightarrow \infty }\frac{t\left( 1-G(\varepsilon \varphi (t))\right)
}{\mathbf{P}\left( \mathcal{H}(t)\right) }=0  \label{Tail2}
\end{equation}%
for any $\varepsilon >0,$ then for any fixed $j\geq 1$ and $y>0$
\begin{equation}
\lim_{t\rightarrow \infty }\mathbf{P}(Z(t-y\varphi (t),t)=j|\mathcal{H}(t))=%
\frac{y}{(j-1)!}\int_{0}^{\frac{1}{y}}z^{j-1}e^{-z}dz.
\end{equation}
\end{theorem}

\textbf{Remark 1}. For the case of the ordinary Galton-Watson processes this
statement was proved in \cite{LV2018}.

\textbf{Remark 2}. It will be shown in Lemma \ref{L_local} below that
\begin{equation*}
\mathbf{P}\left( \mathcal{H}(t)\right) \sim \frac{\varphi (t)}{Bt^{2}}.
\end{equation*}%
Hence (\ref{Tail2}) may be rewritten as

\begin{equation*}
\lim_{t\rightarrow \infty }\frac{t^{3}\left( 1-G(\varepsilon \varphi
(t))\right) }{\varphi (t)}=0.
\end{equation*}

Let
\begin{equation*}
\beta (t):=\max \left\{ 0\leq s<t:Z(s,t)=1\right\}
\end{equation*}%
be the birth moment of the so-called most recent common ancestor (MRCA) of
all particles existing in the population at moment $t$ and let $%
d(t):=t-\beta (t)$ be the distance from the point of observation $t$ to the
birth moment of the MRCA.

Taking $j=1$ in Theorem \ref{T_main} and observing that $\left\{ d(t)\leq
y\varphi (t)\right\} =\left\{ Z(t-y\varphi (t),t)=j\right\} $ we obtain the
following statement.

\begin{corollary}
\label{C_1}If the conditions of Theorem \ref{T_main} are valid then for any $%
y>0$
\begin{equation}
\lim_{t\rightarrow \infty }\mathbf{P}\left( d(t)\leq y\varphi (t)|\mathcal{H}%
(t)\right) =y\left( 1-\exp \left( -\frac{1}{y}\right) \right) .
\end{equation}
\end{corollary}

Our next theorem deals with the case $\varphi (t)=Bat$ for some $a>0.$ Here
much stronger statement may be proved.

\begin{theorem}
\label{T_main2}If Condition A1 is valid, the function $G(t)$ is non-lattice
and
\begin{equation}
\lim_{t\rightarrow \infty }t^{2}\left( 1-G(t)\right) =0,
\end{equation}%
then for any fixed $a>0,$ $j\geq 1$ and $x\in (0,1)$
\begin{equation}
\lim_{t\rightarrow \infty }\mathbf{P}\left( Z(xt,t)=j|0<Z(t)<Bat\right) =%
\frac{1}{(j-1)!}\int_{0}^{\frac{a}{1-x}}z^{j-1}e^{-z}dz\times \frac{\left(
1-x\right) x^{j-1}}{1-e^{-a}}.
\end{equation}
\end{theorem}

Taking $j=1$ in Theorem \ref{T_main2} and observing that $\left\{ d(t)\leq
xt\right\} =\left\{ Z((1-x)t,t)=1\right\} $ we obtain the following
statement:

\begin{corollary}
\label{C_2}If the conditions of Theorem \ref{T_main2} are valid then for any
$x\in (0,1)$
\begin{equation}
\lim_{t\rightarrow \infty }\mathbf{P}\left( d(t)\leq xt|0<Z(t)<Bat\right) =x%
\frac{1-e^{-a/x}}{1-e^{-a}}.
\end{equation}
\end{corollary}

The remaining part of the paper looks as follows. In Section 2 we prove some
auxiliary results. Sections 3 and 4 contains proofs of Theorems \ref{T_main}
and \ref{T_main2}, respectively. We note that Lemmas \ref{L_difference}, \ref%
{L_derivative} and Theorem \ref{T_main2} are proved by V.Vatutin, all other
results are established by Wenming Hong and Yao Ji.

\section{Auxiliary results}

We write
\begin{equation*}
\mathbf{P}(Z(t-By\varphi (t),t)=j|\mathcal{H}(t))=\frac{\mathbf{P}\left(
\mathcal{H}(t)|Z(t-By\varphi (t),t)=j\right) \mathbf{P}\left( Z(t-By\varphi
(t),t)=j\right) }{\mathbf{P}\left( \mathcal{H}(t)\right) }.
\end{equation*}%
Our aim is to investigate separately the asymptotic behavior of each
probability at the right-hand side of this equality.

We start our arguments by the following lemma due to Topchii \cite{Top82}.

\begin{lemma}
\label{L_topchii}If
\begin{equation}
\mathbf{E}\tau ^{3}<\infty ,\,\mathbf{E}\xi =1,\,\sigma ^{2}>0,\,\mathbf{E}%
\xi ^{2}\log (\xi +1)<\infty ,  \label{topLoc}
\end{equation}%
and $G$ is a nondegenerate lattice distribution with span 1 then, as $%
t\rightarrow \infty $
\begin{equation}
t^{2}e^{\frac{k}{Bt}}\mathbf{P}(Z(t)=k)-\frac{1}{B^{2}}\rightarrow 0
\label{LLocal1}
\end{equation}%
uniformly in $0<k\leq Ct<\infty $. Besides, there exists a constant $%
C_{1}<\infty $ such that
\end{lemma}

\begin{equation}
\sup_{k>0,t\geq 0}t^{2}\mathbf{P}(Z(t)=k)\leq C_{1}<\infty .
\label{LocalBound}
\end{equation}%
\hfill \rule{2mm}{3mm}\vspace{4mm}

Note that the condition $\mathbf{E}\tau ^{3}<\infty $ in the lemma cannot be
reduced to $\mu =\mathbf{E}\tau <\infty $. Indeed, if, for instance,%
\begin{equation*}
1-G(t)\sim \frac{c}{t^{\beta }}
\end{equation*}%
as $t\rightarrow \infty $ then, for $1<\beta \leq 2$ and each fixed $k$
there exists%
\begin{equation*}
\lim_{t\rightarrow \infty }t^{\beta /2}\mathbf{P}(Z(t)=k)\in \left( 0,\infty
\right)
\end{equation*}%
(see \cite{Vat1977}), while if $2<\beta <3$ then
\begin{equation*}
\lim_{t\rightarrow \infty }t^{\beta -1}\mathbf{P}(Z(t)=k)\in \left( 0,\infty
\right)
\end{equation*}%
if $k(\beta -1)\leq 1$ and
\begin{equation*}
\lim_{t\rightarrow \infty }t^{2}\mathbf{P}(Z(t)=k)\in \left( 0,\infty \right)
\end{equation*}%
if $k(\beta -1)>1$ (see \cite{Vat1981}).

In what follows we agree to understand (if otherwise is not stated) the
symbol $\sim $\ as $\overset{t\rightarrow \infty }{\sim }$.

\begin{lemma}
\label{L_local}If conditions (\ref{topLoc}) are valid and $G$ is a
nondegenerate lattice distribution with span~1 and $\varphi (t)=o(t)$ as $%
t\rightarrow \infty $ then
\begin{equation}
\mathbf{P}\left( \mathcal{H}(t)\right) \sim \frac{\varphi (t)}{Bt^{2}};
\label{LLocal2}
\end{equation}

2) if the conditions of Theorem \ref{T_main2} are valid then
\begin{equation}
\mathbf{P}\left( 0<Z(t)<Bat\right) \sim \left( 1-e^{-a}\right) \mathbf{P}%
\left( Z(t)>0\right) \sim \frac{1-e^{-a}}{Bt}  \label{LLocal3}
\end{equation}%
for any $a>0.$
\end{lemma}

\textbf{Proof}. Using Lemma \ref{L_topchii} we conclude that%
\begin{equation*}
\mathbf{P}\left( \mathcal{H}(t)|Z(0)=1\right) =\sum_{1\leq k\leq B\varphi
(t)}\mathbf{P}\left( Z(t)=k|Z(0)=1\right) \sim \frac{1}{B^{2}t^{2}}%
\sum_{1\leq k\leq B\varphi (t)}1\sim \frac{\varphi (t)}{Bt^{2}}
\end{equation*}%
proving (\ref{LLocal2}).

To check (\ref{LLocal3}) we recall that by (\ref{Yaglom0})
\begin{equation*}
\lim_{t\rightarrow \infty }\mathbf{P}\left( 0<Z(t)<Bat|Z(t)>0\right)
=1-e^{-a}
\end{equation*}%
and use (\ref{SurvivalProbab1}). \hfill \rule{2mm}{3mm}\vspace{4mm}

Using Lemma \ref{L_topchii} we prove the following statement.

\begin{lemma}
\label{L_difference}Assume that the conditions of Theorem \ref{T_main} are
valid and $\psi (t)\rightarrow \infty $ as $t\rightarrow \infty $ in such a
way that $\psi (t)t^{-1}\rightarrow 0.$ Then
\begin{equation*}
F\left( t;1-\frac{1}{\psi (t)}\right) -F(t;0)\sim \frac{\psi (t)}{B^{2}t^{2}}%
.
\end{equation*}
\end{lemma}

\textbf{Proof}. We have%
\begin{equation*}
F\left( t;1-\frac{1}{\psi (t)}\right) =\sum_{k=0}^{\infty }\,\mathbf{P}%
(Z(t)=k)\left( 1-\frac{1}{\psi (t)}\right) ^{k}.
\end{equation*}%
By the inequality $1-x\leq e^{-x},x\geq 0,$ we conclude that
\begin{equation}
z^{k}=\left( 1-\frac{1}{\psi (t)}\right) ^{k}\leq \exp \left( -\frac{k}{\psi
(t)}\right) .  \notag
\end{equation}%
This and (\ref{LocalBound}) imply for any fixed $N$ and sufficiently large $%
t:$
\begin{align}
\sum_{k>\psi (t)N}\mathbf{P}(Z(t)& =k)z^{k}\leq C_{1}\sum_{k>\psi (t)N}\frac{%
1}{t^{2}}\exp \left( -\frac{k}{\psi (t)}\right)  \notag \\
& =e^{-N}\frac{C_{1}}{t^{2}}\left( 1-e^{-1/\psi (t)}\right) ^{-1}\leq
2e^{-N}C_{1}\frac{\psi (t)}{t^{2}}  \label{Term1}
\end{align}%
and
\begin{equation}
\sum_{0<k<\varepsilon \psi (t)}\mathbf{P}(Z(t)=k)z^{k}\leq
\sum_{0<k<\varepsilon \psi (t)}\mathbf{P}(Z(t)=k)\leq \varepsilon C_{1}\frac{%
\psi (t)}{t^{2}}.  \label{Term2}
\end{equation}

The intermediate term with $\varepsilon \psi (t)<k<N\psi (t)$ is evaluated
as
\begin{align}
\sum_{\varepsilon \psi (t)<k<N\psi (t)}\mathbf{P}(Z(t)& =k)z^{k}\sim \frac{1%
}{(Bt)^{2}}\sum_{\varepsilon \psi (t)<k<N\psi (t)}z^{k}  \notag \\
& \leq \frac{\psi (t)}{(Bt)^{2}}\left( \left( 1-\frac{1}{\psi (t)}\right)
^{\varepsilon \psi (t)}-\left( 1-\frac{1}{\psi (t)}\right) ^{N\psi
(t)+1}\right)  \notag \\
& \sim \frac{\psi (t)}{(Bt)^{2}}\left( e^{-\varepsilon }-e^{-N}\right) .
\label{Term3}
\end{align}%
Combining (\ref{Term1}) - (\ref{Term3}) and letting $\varepsilon \downarrow
0 $ and $N\uparrow \infty $ we obtain the statement of the lemma. \rule%
{2mm}{3mm}\vspace{4mm}

For convenience of references we recall Fa\`{a} di Bruno's formula for the
derivatives of composite functions:

If $i_{r}\in \mathbb{N}_{0}:=\mathbb{N}\cup \left\{ 0\right\} ,r=1,2,...,k$,
$I_{k}:=i_{1}+\cdots +i_{k}$ and
\begin{equation*}
\mathcal{D}(k):=\left\{ \left( i_{1},...,i_{k}\right) :1\cdot i_{1}+2\cdot
i_{2}+\cdot \cdot \cdot +ki_{k}=k\right\} ,
\end{equation*}%
then for the derivatives of the composition $H(T(z))$ of the functions $%
H(\cdot )$ and $T(\cdot )$ we have%
\begin{equation}
\frac{d^{k}}{dz^{k}}\left[ H(T(z))\right] =\sum_{\mathcal{D}(k)}\frac{k!}{%
i_{1}!\cdot \cdot \cdot i_{k}!}H^{(I_{k})}(T(z))\prod_{r=1}^{k}\left( \frac{%
T^{(r)}(z)}{r!}\right) ^{i_{r}}.  \label{Faa}
\end{equation}%
The next lemma is crucial for the proof of Theorem \ref{T_main}.

\begin{lemma}
\label{L_derivative} If $\psi (t)\rightarrow \infty $ as $t\rightarrow
\infty $ in such a way that $\psi (t)t^{-1}\rightarrow 0$ and the conditions
of Theorem \ref{T_main} are valid then,m for any fixed $k\in \mathbb{N}$
\begin{equation}
F^{(k)}(t;f(F(\psi (t))))\sim \frac{(B\psi (t))^{k+1}}{B^{2}t^{2}}k!.
\label{AsymDer}
\end{equation}
\end{lemma}

\textbf{Proof}. It follows from Lemma \ref{L_difference} that, for any
positive $\lambda $%
\begin{equation}
\lim_{t\rightarrow \infty }\frac{t^{2}}{\psi (t)}\left[ F\left( t;1-\frac{%
\lambda }{B\psi (t)}\right) -F(t;0)\right] =\frac{1}{B\lambda }.  \notag
\end{equation}%
Set for brevity $F(t):=F(t;0)$ and take a fixed $\lambda >0.$ Since $%
f^{\prime }(1)=1$ and $f(F(\psi (t)))\rightarrow 1$ as $t\rightarrow \infty
, $ it follows that
\begin{eqnarray*}
1-f^{\lambda }(F(\psi (t))) &=&1-(1-\left( 1-f(F(\psi (t)))\right) ^{\lambda
} \\
&\sim &\lambda \left( 1-f(F(\psi (t)))\right) \sim \lambda \left( 1-F(\psi
(t))\right) \sim \frac{\lambda }{B\psi (t)}.
\end{eqnarray*}%
Hence, setting for brevity $w(t):=f(F(\psi (t)))$ we get for any positive $%
\lambda $

\begin{equation}
\lim_{t\rightarrow \infty }\frac{t^{2}}{\psi (t)}\left[ F(t;w^{\lambda
}(t))-F(t;0)\right] =\frac{1}{B\lambda }.  \label{Analyt1}
\end{equation}%
Since the prelimiting and limiting functions in (\ref{Analyt1}) \ are
analytical in the complex domain $Re$ $\lambda >0,$ the derivatives of any
order of the prelimiting functions with respect to $\lambda $ converge to
the derivatives of the respective order of the limiting function. Hence it
follows that%
\begin{eqnarray}
&&\lim_{t\rightarrow \infty }\frac{\partial ^{k}}{\partial \lambda ^{k}}%
\left( \frac{t^{2}}{\psi (t)}\left[ F\left( t;w^{\lambda }(t)\right) -F(t;0)%
\right] \right) =\lim_{t\rightarrow \infty }\frac{t^{2}}{\psi (t)}\frac{%
\partial ^{k}}{\partial \lambda ^{k}}F\left( t;w^{\lambda }(t)\right)  \notag
\\
&&\qquad \qquad \qquad \qquad \qquad \qquad \quad =\lim_{t\rightarrow \infty
}\frac{t^{2}}{\psi (t)}\frac{\partial ^{k}}{\partial \lambda ^{k}}F\left(
t;w^{\lambda }(t)\right) =\left( -1\right) ^{k}\frac{k!}{B\lambda ^{k+1}}.
\label{Derivat3}
\end{eqnarray}%
In particular, for $k=1$%
\begin{eqnarray}
\frac{t^{2}}{\psi (t)}\frac{\partial }{\partial \lambda }F\left(
t;w^{\lambda }(t)\right) &=&\frac{t^{2}}{\psi (t)}F^{\prime }\left(
t;w^{\lambda }(t)\right) w^{\lambda }(t)\log w(t)  \notag \\
&\sim &\left( -1\right) \frac{1}{B\lambda ^{2}}.  \label{FirstDerivat}
\end{eqnarray}%
Using the equivalences $\log (1-x)\sim -x$ as $x\downarrow 0$ and
\begin{equation}
1-f\left( F(\psi (t))\right) \sim 1-F(\psi (t))=\mathbf{P}\left( Z(\psi
(t))>0\right) \sim \frac{1}{B\psi (t)}  \label{Equival2}
\end{equation}%
as $t\rightarrow \infty $ it is not difficult to deduce from (\ref%
{FirstDerivat}) with $\lambda =1$ that%
\begin{equation*}
F^{\prime }\left( t;w(t)\right) \sim \frac{\psi ^{2}(t)}{t^{2}}=\frac{\left(
B\psi (t)\right) ^{2}}{B^{2}t^{2}}\times 1!
\end{equation*}%
proving the lemma for $k=1$.

Assume that the asymptotic representation
\begin{equation*}
F^{(r)}(t;w(t))\sim \frac{(B\psi (t))^{r+1}}{B^{2}t^{2}}r!
\end{equation*}%
is valid for all $r<k$. By Fa\`{a} di Bruno's formula (\ref{Faa}) we have%
\begin{eqnarray*}
&&\frac{t^{2}}{\psi (t)}\frac{\partial ^{k}}{\partial \lambda ^{k}}%
[F(t;w^{\lambda }(t)] \\
&=&\frac{t^{2}}{\psi (t)}\sum_{\mathcal{D}(k)}\frac{k!}{i_{1}!\cdots i_{k}!}%
F^{(I_{k})}(t;w^{\lambda }(t))\prod_{r=1}^{k}\left( \frac{1}{r!}\frac{%
\partial ^{r}}{\partial \lambda ^{r}}w^{\lambda }(t)\right) ^{i_{r}} \\
&=&\frac{t^{2}}{\psi (t)}\sum_{\mathcal{D}(k)}\frac{k!}{i_{1}!\cdots i_{k}!}%
F^{(I_{k})}(t;w^{\lambda }(t))\prod_{r=1}^{k}\left( \frac{w^{\lambda }(t)}{r!%
}\log ^{r}w(t)\right) ^{i_{r}} \\
&=&\frac{t^{2}}{\psi (t)}\log ^{k}w(t)\sum_{\mathcal{D}(k)}\frac{k!}{%
i_{1}!\cdots i_{k}!}F^{(I_{k})}(t;w^{\lambda }(t))w^{\lambda
I_{k}}(t)\prod_{r=1}^{k}\frac{1}{(r!)^{i_{r}}} \\
&\sim &(-1)^{k}\frac{t^{2}}{\psi (t)}\left( \frac{1}{B\psi (t)}\right)
^{k}\sum_{\mathcal{D}(k)}\frac{k!}{i_{1}!\cdots i_{k}!}F^{(I_{k})}(t;w^{%
\lambda }(t))w^{\lambda I_{k}}(t)\prod_{r=1}^{k}\frac{1}{(r!)^{i_{r}}}.
\end{eqnarray*}%
Set $\mathcal{D}^{\prime }(k)=\mathcal{D}(k)\backslash \left\{
(k,...,0,0)\right\} .$ In view of induction hypothesis (recall (\ref{AsymDer}%
)) and the estimate $I_{k}=i_{1}+\cdots +i_{k}\leq k-1$ valid for all \ $%
\left( i_{1},...,i_{k}\right) \in \mathcal{D}^{\prime }(k)$ we see that
\begin{eqnarray*}
&&\lim_{t\rightarrow \infty }\frac{t^{2}}{\psi (t)}\left( \frac{1}{B\psi (t)}%
\right) ^{k}\sum_{\mathcal{D}^{\prime }(k)}\frac{k!}{i_{1}!\cdots i_{k}!}%
F^{(I_{k})}(t;f(F(\psi (t))))f^{I_{k}}(F(\psi (t)))\prod_{r=1}^{k}\frac{1}{%
(r!)^{i_{r}}} \\
&&\qquad \qquad \qquad \quad \leq \sum_{\mathcal{D}^{\prime }(k)}\frac{k!}{%
i_{1}!\cdots i_{k}!}\lim_{t\rightarrow \infty }\frac{t^{2}}{\psi (t)}\left(
\frac{1}{B\psi (t)}\right) ^{k}F^{(I_{k})}(t;f(F(\psi (t))))=0.
\end{eqnarray*}%
Hence, setting $\lambda =1$ we conclude by (\ref{Derivat3}) that%
\begin{eqnarray*}
\left( -1\right) ^{k}\frac{k!}{B} &\sim &\frac{t^{2}}{\psi (t)}\frac{%
\partial ^{k}}{\partial \lambda ^{k}}[F(t;w^{\lambda }(t))]\Big|_{\lambda =1}
\\
&=&o\left( 1\right) +\frac{t^{2}}{\psi (t)}\left( \frac{1}{B\psi (t)}\right)
^{k}(-1)^{k}\frac{k!}{k!\cdots 0!0!}F^{(k)}(t;w(t))w^{k}(t) \\
&=&o\left( 1\right) +\frac{t^{2}}{\psi (t)}\left( \frac{1}{B\psi (t)}\right)
^{k}(-1)^{k}F^{(k)}(t;w(t)).
\end{eqnarray*}%
Therefore,%
\begin{equation*}
F^{(k)}(t;f(F(\psi (t))))=F^{(k)}(t;w(t))\sim \frac{(B\psi (t))^{k+1}}{%
B^{2}t^{2}}k!
\end{equation*}%
that completes the induction step and proves Lemma \ref{L_derivative}.

\vspace{4mm} We now consider a Bellman-Harris branching process which is
initiated at time $t=0$ by a \textit{random} number of particles distributed
the same as $\xi $ specified by $f(s)$ in (\ref{DefGener}). The initial
particles as well as the other particles have life-length distribution $G(t)$%
. Each particle of the process produces children at the end of its life in
accordance with probability generating function $f(s)$. We denote this new
process as $Y(t)$. Clearly,%
\begin{equation*}
T(t;s):=\mathbf{E}\left[ s^{Y(t)}\right] =f\left( F(t;s)\right)
\end{equation*}%
and, as a result%
\begin{equation*}
\mathbf{P}\left( Y(t)>0\right) =1-f\left( F(t;0)\right) \sim 1-F(t;0)\sim
\frac{1}{Bt}
\end{equation*}%
and, in view of%
\begin{eqnarray*}
\mathbf{E}\left[ e^{-2\mu \lambda Y(t)/\sigma ^{2}t}|Y(t)>0\right] &=&\frac{%
\mathbf{E}\left[ e^{-2\mu \lambda Y(t)/\sigma ^{2}t};Y(t)>0\right] }{\mathbf{P%
}\left( Y(t)>0\right) } \\
&=&\frac{T\left( t;e^{-2\mu \lambda /\sigma ^{2}t}\right) -T(t;0)}{%
1-T(t;0)} \\
&=&\frac{f\left( F\left( t;e^{-2\mu \lambda /\sigma ^{2}t}\right) \right)
-f\left( F(t;0)\right) }{1-f\left( F(t;0)\right) }
\end{eqnarray*}%
and (\ref{Yaglom0})
\begin{eqnarray}
\lim_{t\rightarrow \infty }\mathbf{E}\left[ e^{-2\mu \lambda Y(t)/\sigma
^{2}t}|Y(t)>0\right] &=&1-\lim_{t\rightarrow \infty }\frac{1-f\left( F\left(
t;e^{-2\mu \lambda /\sigma ^{2}t}\right) \right) }{1-f\left(
F(t;0)\right) }  \notag \\
&=&1-\lim_{t\rightarrow \infty }\frac{1-F\left( t;e^{-2\mu \lambda
/\sigma ^{2}t}\right) }{1-F(t;0)}  \notag \\
&=&\lim_{t\rightarrow \infty }\mathbf{E}\left[ e^{-2\mu \lambda Z(t)/\sigma
^{2}t}|Z(t)>0\right] =\frac{1}{1+\lambda }.  \label{YaglomExp2}
\end{eqnarray}%
Hence, the limiting conditional distribution of the process $Y(t)$ given $%
\left\{ Y(t)>0\right\} $ is exponential with parameter 1.\ \ \

Let $Z^{\ast }(t,x)$ be the number of particles existing in the process at
moment $t$, which will exist at moment $t+x$.

The following statement, showing that under the conditions of Theorem \ref%
{T_main2} the probability that there is a particle at time $t$ which will
survive up to moment $t+\varepsilon t$ is negligible with $\mathbf{P}%
(Z(t)>0),$ is a particular case of Lemma 1 in \cite{Vat1979}.

\begin{lemma}
\label{L_neglig2}If the conditions of Theorem \ref{T_main2} are valid then
for any $\varepsilon >0$,
\begin{equation}
\lim_{t\rightarrow \infty }\frac{\mathbf{P}(Z^{\ast }(t,\varepsilon t)>0)}{%
\mathbf{P}(Z(t)>0)}=0.
\end{equation}
\end{lemma}

We complement Lemma \ref{L_neglig2} by the following result:

\begin{lemma}
\label{L_neglig}If the conditions of Theorem \ref{T_main} are valid then for
any $\varepsilon >0$,
\begin{equation}
\lim_{t\rightarrow \infty }\frac{\mathbf{P}(Z^{\ast }(t,\varepsilon \varphi
(t))>0)}{\mathbf{P}(\mathcal{H}(t))}=0.
\end{equation}
\end{lemma}

\textbf{Proof}. Let $\tilde{Z}(t,x)$ be the number of particles at moment $t$
whose age does not exceed $x$. Setting
\begin{equation*}
F(t,x;s):=\mathbf{E}\left[ s^{\tilde{Z}(t,x)}|Z(0)=1\right]
\end{equation*}%
and introducing the notation $J(y)=1$ for $y\geq 0$ and $J(y)=0$ for $y<0,$
we deduce by the total probability formula the integral equation
\begin{equation}
F(t,x;s)=(1-G(t))[sJ(x-t)+1-J(x-t)]+\int_{0}^{t}f(F(t-u,x;s))dG(u).  \notag
\end{equation}

Denoting $A(t,x):=\mathbf{E}\tilde{Z}(t,x)$ we conclude by the previous
relation that
\begin{equation}
A(t,x)=(1-G(t))J(x-t)+\int_{0}^{t}A(t-u,x)dG(u).  \notag
\end{equation}%
Solving this renewal type equation gives
\begin{equation}
A(t,x)=\int_{0}^{t}(1-G(t-u))J(x-(t-u))d{U}(u),  \notag
\end{equation}%
where $U(t)=\sum\limits_{k=0}^{\infty }G^{\ast k}(t).$ In particular,

\begin{equation}
\mathbf{E}Z(t)=A(t)=A(t,t)=1=\int_{0}^{t}(1-G(t-u))dU(u).  \notag
\end{equation}

We know that
\begin{equation*}
\mathbf{E}[Z(t+\varepsilon \varphi (t))]=1=\int_{0}^{t+\varepsilon \varphi
(t)}(1-G(t+\varepsilon \varphi (t)-u)dU(u)
\end{equation*}%
and
\begin{align*}
\mathbf{E}[\tilde{Z}(t+\varepsilon \varphi (t),\varepsilon \varphi (t))]&
=A(t+\varepsilon \varphi (t),\varepsilon \varphi (t)) \\
& =\int_{0}^{t+\varepsilon \varphi (t)}(1-G(t+\varepsilon \varphi
(t)-u)J(\varepsilon \varphi (t)-(t+\varphi (t)\varepsilon )-u)dU(u) \\
& =\int_{t}^{t+\varepsilon \varphi (t)}(1-G(t+\varepsilon \varphi
(t)-u))dU(u).
\end{align*}

Since

\begin{equation}
Z^{\ast }(t,\varepsilon \varphi (t))=Z(t+\varepsilon \varphi (t))-\tilde{Z}%
(t+\varepsilon \varphi (t),\varepsilon \varphi (t))  \notag
\end{equation}%
for any $\varepsilon >0,$ it follows by Markov inequality that

\begin{align*}
\mathbf{P}(Z^{\ast }(t,\varepsilon \varphi (t)& \geq 1)\leq \mathbf{E}%
Z^{\ast }(t,\varepsilon \varphi (t))=\mathbf{E}[Z(t+\varepsilon \varphi (t))-%
\tilde{Z}(t+\varepsilon \varphi (t),\varepsilon \varphi (t))] \\
& =\int_{0}^{t}(1-G(t+\varepsilon \varphi (t)-u))dU(u) \\
& \leq U(t)(1-G(\varepsilon \varphi (t)))\leq C\frac{t}{\mu }\left(
1-G(\varepsilon \varphi (t))\right) =o\left( \mathbf{P}(\mathcal{H}%
(t))\right)
\end{align*}%
in view of (\ref{Tail2}) and the asymptotic relation $U(t)\sim t\mu ^{-1}$
as $t\rightarrow \infty $ being valid by the key renewal theorem for the
renewal function $U(t)$ with finite mean $\mu $ for the increments.

Lemma \ref{L_neglig} is proved.\

\section{Proof of Theorem \protect\ref{T_main}}

Let $\zeta _{i}:=\zeta _{i}(t-y\varphi (t)),\,i=1,2,...,Z(t-y\varphi (t))$
be the remaining life-lengths of the particles existing in the process at
moment $t-y\varphi (t)$. We fix $\varepsilon >0$ and $y>0$ and introduce the
event
\begin{equation*}
\mathcal{C}\left( t,y,\varepsilon \right) :=\left\{ \max_{1\leq i\leq
Z(t-y\varphi (t))}\zeta _{i}\leq \varepsilon \varphi (t)\right\}
\end{equation*}%
and the event $\mathcal{\bar{C}}\left( t,y,\varepsilon \right) $
complementary to $\mathcal{C}\left( t,y,\varepsilon \right) $. In view of 
Lemma \ref{L_neglig} and monotonicity of $\varphi (t)$ %
\begin{equation}
\lim_{t\rightarrow \infty }\frac{\mathbf{P}\left( \mathcal{\bar{C}}\left(
t,y,\varepsilon \right) \right) }{\mathbf{P}\left( \mathcal{H}(t)\right) }
=\lim_{t\rightarrow \infty }\frac{\mathbf{P}\left( \mathcal{\bar{C}}\left(
t,y,\varepsilon \right) \right) }{\mathbf{P}\left( \mathcal{H}(t-y\varphi (t))\right) }\frac{\mathbf{P}\left( \mathcal{H}(t-y\varphi (t))\right) }{\mathbf{P}\left( \mathcal{H}(t)\right)}=0.
\label{NegligC}
\end{equation}
Thus, for any $j\geq 1$
\begin{equation}
\mathbf{P}(Z(t-y\varphi (t),t)=j)=\mathbf{P}(Z(t-y\varphi (t),t)=j;\mathcal{C%
}\left( t,y,\varepsilon \right) )+o\left( \mathbf{P}\left( \mathcal{H}%
(t)\right) \right) .  \label{Reduced2}
\end{equation}%
Set%
\begin{equation*}
\mathcal{C}_{k}\left( t,y,\varepsilon \right) :=\mathcal{C}\left(
t,y,\varepsilon \right) \cap \left\{ Z(t-y\varphi (t))=k\right\} .
\end{equation*}%
Then, for $k\geq j$
\begin{align*}
& \mathbf{P}(\mathcal{C}_{k}\left( t,y,\varepsilon \right) ;Z(t-y\varphi
(t),t)=j) \\
& =\mathbf{E}\left[ \mathbf{P}(\mathcal{C}_{k}\left( t,y,\varepsilon \right)
;Z(t-y\varphi (t),t)=j\ |\ \zeta _{i},i=1,2,...,k)\right] \\
& =\mathbf{P}(\mathcal{C}_{k}\left( t,y,\varepsilon \right) ) \\
\times & \mathbf{E}[\sum_{0\leq i_{1}<i_{2}<...<i_{j}\leq k}\prod_{i\in
\left\{ i_{1},...i_{j}\right\} }(1-f(F(y\varphi (t)-\zeta
_{i})))\prod_{1\leq i\leq k:i\notin \left\{ i_{1},...i_{j}\right\}
}f(F(y\varphi (t)-\zeta _{i}))|\mathcal{C}_{k}\left( t,y,\varepsilon \right)
] \\
& \qquad \geqslant \mathbf{P}(\mathcal{C}_{k}\left( t,y,\varepsilon \right)
)C_{k}^{j}\mathbf{E}\left[ (1-f(F(y\varphi (t))))^{j}f^{k-j}(F(\left(
y-\varepsilon \right) \varphi (t)))|\mathcal{C}_{k}\left( t,y,\varepsilon
\right) \right] \\
& \qquad =\mathbf{P}(\mathcal{C}_{k}\left( t,y,\varepsilon \right)
)C_{k}^{j}(1-f(F(y\varphi (t))))^{j}f^{k-j}(F(\left( y-\varepsilon \right)
\varphi (t))) \\
& \qquad \geq \mathbf{P}(Z(t-y\varphi (t))=k)C_{k}^{j}(1-f(F(y\varphi
(t))))^{j}f^{k-j}(F(\left( y-\varepsilon \right) \varphi (t))) \\
& \qquad \qquad -\mathbf{P}(Z(t-y\varphi (t)=k,\mathcal{\bar{C}}\left(
t,y,\varepsilon \right) ).
\end{align*}

By the same arguments we get,
\begin{eqnarray*}
\mathbf{P}(\mathcal{C}_{k}\left( t,y,\varepsilon \right) ,Z(t-y\varphi
(t),t)=j) &\leq &\mathbf{P}(\mathcal{C}_{k}\left( t,y,\varepsilon \right)
)C_{k}^{j}(1-f(F(\left( y-\varepsilon \right) \varphi
(t))))^{j}f^{k-j}(F(y\varphi (t))) \\
&\leq &\mathbf{P}(Z(t-y\varphi (t))=k)C_{k}^{j}(1-f(F(\left( y-\varepsilon
\right) \varphi (t))))^{j}f^{k-j}(F(y\varphi (t))).
\end{eqnarray*}

As a result we obtain
\begin{align}
\mathbf{P}(Z(t-y\varphi (t),t)=j,\mathcal{C}\left( t,y,\varepsilon \right)
)& =\sum_{k=j}^{\infty }\mathbf{P}(\mathcal{C}_{k}\left( t,y,\varepsilon
\right) ;Z(t-y\varphi (t),t)=j)  \notag \\
& \leq \sum_{k=j}^{\infty }\mathbf{P}(Z(t-y\varphi
(t))=k)C_{k}^{j}(1-f(F(\left( y-\varepsilon \right) \varphi
(t))))^{j}f^{k-j}(F(y\varphi (t)))  \notag \\
& =\frac{(1-f(F(\left( y-\varepsilon \right) \varphi (t)))))^{j}}{j!}%
F^{(j)}(t-y\varphi (t);f(F(y\varphi (t))))  \label{SemiAbove}
\end{align}%
and
\begin{align}
\mathbf{P}(Z(t-y\varphi (t),t)=j,\mathcal{C}\left( t,y,\varepsilon \right)
)& \geq \sum_{k=j}^{\infty }C_{k}^{j}(1-f(F(y\varphi
(t))))^{j}f^{k-j}(F(\left( y-\varepsilon \right) \varphi (t)))\mathbf{P}%
(Z(t-y\varphi (t))=k)  \notag \\
& \qquad \qquad \qquad \qquad \qquad -\sum_{k=j}^{\infty }\mathbf{P}%
(Z(t-y\varphi (t))=k,\mathcal{\bar{C}}\left( t,y,\varepsilon \right) )
\notag \\
& \geq \frac{(1-f(F(y\varphi (t))))^{j}}{j!}F^{(j)}(t-y\varphi
(t);f(F(\left( y-\varepsilon \right) \varphi (t))))-\mathbf{P}(\mathcal{\bar{%
C}}\left( t,y,\varepsilon \right) ).  \label{SemiBelow}
\end{align}%
We know by Lemma \ref{L_derivative} that
\begin{equation*}
F^{(j)}(t-y\varphi (t);f(F(y\varphi (t))))\sim \frac{(By\varphi (t))^{j+1}}{%
B^{2}t^{2}}j!.
\end{equation*}%
Thus, in view of (\ref{SemiAbove})
\begin{align}
& \limsup_{t\rightarrow \infty }\frac{\mathbf{P}(Z(t-y\varphi (t),t)=j,%
\mathcal{C}\left( t,y,\varepsilon \right) )}{\mathbf{P}(\mathcal{H}(t))}
\notag \\
& \qquad \leq \limsup_{t\rightarrow \infty }\frac{Bt^{2}}{\varphi (t)}\times
\frac{1}{B}y\varphi (t)\left( \frac{y\varphi (t)}{\varphi (t)(y-\varepsilon )%
}\right) ^{j}\frac{1}{t^{2}}=y\left( \frac{y}{y-\varepsilon }\right) ^{j}
\notag
\end{align}%
and by (\ref{SemiBelow})
\begin{align*}
& \liminf_{t\rightarrow \infty }\frac{\mathbf{P}(Z(t-y\varphi (t),t)=j,%
\mathcal{C}\left( t,y,\varepsilon \right) )}{\mathbf{P}(\mathcal{H}(t))} \\
& \qquad \geq \liminf_{t\rightarrow \infty }\left[ \frac{t^{2}B}{\varphi (t)}%
\times \frac{1}{B}(\varphi (t)(y-\varepsilon ))\left( \frac{\varphi
(t)(y-\varepsilon )}{y\varphi (t)}\right) ^{j}\frac{1}{t^{2}}-\frac{\mathbf{P%
}(\mathcal{\bar{C}}\left( t,y,\varepsilon \right) )}{\mathbf{P}(\mathcal{H}%
(t))}\right] \\
& \qquad \qquad =(y-\varepsilon )\left( \frac{y-\varepsilon }{y}\right) ^{j}.
\end{align*}%
Hence, letting $\varepsilon \rightarrow 0$, we conclude

\begin{equation*}
\lim_{t\rightarrow \infty }\frac{\mathbf{P}(Z(t-y\varphi (t),t)=j)}{\mathbf{P%
}(\mathcal{H}(t))}=\lim_{\varepsilon \rightarrow 0}\lim_{t\rightarrow \infty
}\frac{\mathbf{P}(Z(t-y\varphi (t),t)=j,\mathcal{C}\left( t,y,\varepsilon
\right) )}{\mathbf{P}(\mathcal{H}(t))}=y.
\end{equation*}

Let now $Y_{1}^{\ast }(t),...,Y_{j}^{\ast }(t)$ be a tuple of i.i.d.random
variables distributed as $\{Y(t)|Y(t)>0\}$, and let $\eta _{1},...,\eta _{j}$
be i.i.d.random variables having exponential distributed with parameter 1.
It follows that
\begin{eqnarray*}
\lim_{t\rightarrow \infty }\mathbf{P}(\mathcal{H}(t)|Z(t-y\varphi (t),t)
&=&j;\mathcal{C}\left( t,y,\varepsilon \right) ) \\
&=&\lim_{t\rightarrow \infty }\mathbf{P}\left( \sum_{i=1}^{j}Y_{i}^{\ast
}(y\varphi (t)-\zeta _{i})\leq B\varphi (t)\big|\mathcal{C}\left(
t,y,\varepsilon \right) \right) \\
&=&\lim_{t\rightarrow \infty }\mathbf{P}\left( \sum_{i=1}^{j}\frac{%
Y_{1}^{\ast }(y\varphi (t)-\zeta _{i})}{B\left( y\varphi (t)-\zeta
_{i}\right) }\frac{\left( y\varphi (t)-\zeta _{i}\right) }{y\varphi (t)}\leq
\frac{1}{y}\Big|\mathcal{C}\left( t,y,\varepsilon \right) \right) .
\end{eqnarray*}%
Since%
\begin{equation*}
\mathbf{P}\left( \sum_{i=1}^{j}\frac{Y_{1}^{\ast }(y\varphi (t)-\zeta _{i})}{%
B\left( y\varphi (t)-\zeta _{i}\right) }\frac{\left( y\varphi (t)-\zeta
_{i}\right) }{y\varphi (t)}\leq \frac{1}{y}\Big|\mathcal{C}\left(
t,y,\varepsilon \right) \right) \leq \mathbf{P}\left( \sum_{i=1}^{j}\frac{%
Y_{1}^{\ast }(y\varphi (t)-\zeta _{i})}{B\left( y\varphi (t)-\zeta
_{i}\right) }\frac{y-\varepsilon }{y}\leq \frac{1}{y}\Big|\mathcal{C}\left(
t,y,\varepsilon \right) \right)
\end{equation*}%
and
\begin{equation*}
\mathbf{P}\left( \sum_{i=1}^{j}\frac{Y_{1}^{\ast }(y\varphi (t)-\zeta _{i})}{%
B\left( y\varphi (t)-\zeta _{i}\right) }\frac{\left( y\varphi (t)-\zeta
_{i}\right) }{y\varphi (t)}\leq \frac{1}{y}\Big|\mathcal{C}\left(
t,y,\varepsilon \right) \right) \geq \mathbf{P}\left( \sum_{i=1}^{j}\frac{%
Y_{1}^{\ast }(y\varphi (t)-\zeta _{i})}{B\left( y\varphi (t)-\zeta
_{i}\right) }\leq \frac{1}{y}\Big|\mathcal{C}\left( t,y,\varepsilon \right)
\right) ,
\end{equation*}%
we conclude by (\ref{Yaglom0}) that
\begin{eqnarray*}
&&\lim_{\varepsilon \rightarrow 0}\lim_{t\rightarrow \infty }\mathbf{P}%
\left( \sum_{i=1}^{j}\frac{Y_{1}^{\ast }(y\varphi (t)-\zeta _{i})}{B\left(
y\varphi (t)-\zeta _{i}\right) }\frac{\left( y\varphi (t)-\zeta _{i}\right)
}{y\varphi (t)}\leq \frac{1}{y}\Big|\mathcal{C}\left( t,y,\varepsilon
\right) \right) \\
&&\qquad \qquad \qquad \qquad =\mathbf{P}\left( \sum_{i=1}^{j}\eta _{\iota
}\leq \frac{1}{y}\right) =\frac{1}{(j-1)!}\int_{0}^{\frac{1}{y}%
}z^{j-1}e^{-z}dz.
\end{eqnarray*}

Combining this result with Lemma \ref{L_local} we see that,

\begin{align}
& \lim_{t\rightarrow \infty }\mathbf{P}(Z(t-y\varphi (t),t)=j|\mathcal{H}%
(t))=\lim_{t\rightarrow \infty }\frac{\mathbf{P}(Z(t-y\varphi (t),t)=j)%
\mathbf{P}(\mathcal{H}(t)|Z(t-y\varphi (t),t)=j)}{\mathbf{P}(\mathcal{H}(t))}
\notag \\
& \qquad \qquad \qquad =\lim_{\varepsilon \rightarrow 0}\lim_{t\rightarrow
\infty }\frac{\mathbf{P}(Z(t-y\varphi (t),t)=j,\mathcal{C}\left(
t,y,\varepsilon \right) )\mathbf{P}(\mathcal{H}(t)|Z(t-y\varphi (t),t)=j,%
\mathcal{C}\left( t,y,\varepsilon \right) )}{\mathbf{P}(\mathcal{H}(t))}
\notag \\
& \qquad \qquad \qquad \qquad =\frac{y}{(j-1)!}\int_{0}^{\frac{1}{y}%
}z^{j-1}e^{-z}dz.  \notag
\end{align}%
\hfill \rule{2mm}{3mm}\vspace{4mm}

\section{Proof of Theorem \protect\ref{T_main2}}

The proof of the theorem follows the line of proving Theorem \ref{T_main}.

Let $\zeta _{i}:=\zeta _{i}(xt),\,i=1,2,...,Z(xt)$ be the remaining
life-lengths of the particles existing in the process at moment $t(1-x),x\in
(0,1)$. We fix $\varepsilon >0$ and introduce the event
\begin{equation*}
\mathcal{D}\left( t,x,\varepsilon \right) :=\left\{ \max_{1\leq i\leq
Z(xt)}\zeta _{i}\leq \varepsilon t\right\} .
\end{equation*}
It follows from Lemma 2 in \cite{Sag1986} and Lemma \ref{L_neglig2} of the
present paper that under the conditions of Theorem \ref{T_main2}%
\begin{equation*}
\lim_{t\rightarrow \infty }\frac{\mathbf{P}(Z(xt,t)=j)}{\mathbf{P}(Z(t)>0)}%
=\lim_{\varepsilon \rightarrow 0}\lim_{t\rightarrow \infty }\frac{\mathbf{P}%
(Z(xt,t)=j;\mathcal{D}\left( t,x,\varepsilon \right) )}{\mathbf{P}(Z(t)>0)}%
=\left( 1-x\right) x^{j-1}
\end{equation*}%
for any $j\geq 1$. Now using the arguments similar to those used in the
proof of Theorem \ref{T_main} we have

\begin{eqnarray*}
\lim_{t\rightarrow \infty }\mathbf{P}(0 &<&Z(t)<Bat|Z(xt,t)=j) \\
&=&\lim_{\varepsilon \rightarrow 0}\lim_{t\rightarrow \infty }\mathbf{P}%
(0<Z(t)<Bat|Z(xt,t)=j;\mathcal{D}\left( t,x,\varepsilon \right) ) \\
&=&\lim_{\varepsilon \rightarrow 0}\lim_{t\rightarrow \infty }\mathbf{P}%
\left( \sum_{i=1}^{j}Y_{i}^{\ast }((1-x)t-\zeta _{i})\leq Bat\big|\mathcal{D}%
\left( t,x,\varepsilon \right) \right) \\
&=&\lim_{\varepsilon \rightarrow 0}\lim_{t\rightarrow \infty }\mathbf{P}%
\left( \sum_{i=1}^{j}\frac{Y_{1}^{\ast }((1-x)t-\zeta _{i})}{B\left(
(1-x)t-\zeta _{i}\right) }\frac{\left( (1-x)t-\zeta _{i}\right) }{(1-x)t}%
\leq \frac{a}{1-x}\Big|\mathcal{D}\left( t,x,\varepsilon \right) \right) \\
&=&\mathbf{P}\left( \sum_{i=1}^{j}\eta _{\iota }\leq \frac{a}{1-x}\right) =%
\frac{1}{(j-1)!}\int_{0}^{\frac{a}{1-x}}z^{j-1}e^{-z}dz.
\end{eqnarray*}%
Combining this result with Lemma \ref{L_local} we see that,

\begin{align*}
& \lim_{t\rightarrow \infty }\mathbf{P}(Z(xt,t)=j|0<Z(t)<Bat) \\
& =\lim_{t\rightarrow \infty }\frac{\mathbf{P}(Z(xt,t)=j)\mathbf{P}%
(0<Z(t)<Bat|Z(xt,t)=j)}{\mathbf{P}(0<Z(t)<Bat)} \\
& \qquad \qquad \qquad =\lim_{\varepsilon \rightarrow 0}\lim_{t\rightarrow
\infty }\frac{\mathbf{P}(Z(xt,t)=j,\mathcal{D}\left( t,y,\varepsilon \right)
)\mathbf{P}((0<Z(t)<Bat|Z(xt,t)=j,\mathcal{D}\left( t,y,\varepsilon \right) )%
}{\mathbf{P}(0<Z(t)<Bat)} \\
& \qquad \qquad \qquad \qquad =\frac{1}{(j-1)!}\int_{0}^{\frac{a}{1-x}%
}z^{j-1}e^{-z}dz\times \left( 1-x\right) x^{j-1}\times \left(
1-e^{-a}\right) ^{-1}.
\end{align*}%
Theorem \ref{T_main2} is proved.

To prove Corollary \ref{C_2} we set $j=1$ in the preceding formula and
obtain
\begin{eqnarray*}
\lim_{t\rightarrow \infty }\mathbf{P}(d(t) &\leq
&xt|0<Z(t)<Bat)=\lim_{t\rightarrow \infty }\mathbf{P}%
(Z(t(1-x),t)=j|0<Z(t)<Bat) \\
&=&x\frac{1-e^{-a/x}}{1-e^{-a}}.
\end{eqnarray*}

\end{document}